# Few more Comments on Benford's Law

Jacek  M. Kowalski *

*( Univ. of North Texas, Denton)*

*(retired Assoc.Prof. of  Physics)*

It is pointed out that the language of quotient groups and "wrapped" distributions allows an elementary discussion of Benford's Law, and adds arguments supporting wide-spread observability of this statistics.

There is an vast body of works devoted to Benford's law and its applications- an excellent online bibliography [1]  lists almost one thousand , often collaborating, authors. It is composed and updated by active researchers and authors of an impressive monograph [2] on the subject. Below we add few simple comments on this theme which could have at least some pedagogical value.

For $R^+$ considered as an abelian group under multiplication and given positive integer $b \geq 2$, let $B_b^+$ be its subgroup of integer powers of $b$: $B_b^+ = \{x \in R^+;\ x = b^k.\ k \in Z\}$. Denote by $Q_b$ the quotient group $R^+ / B_b^+$, $Q_b = \{x \in [1,b),\ \bullet\ \mathrm{mod}\, b\}$ , an abelian group under multiplication modulo $b$. Each element $x$ of $Q_b$ represents an infinite sequence of real numbers which elements can be written as $xb^k,\ k \in Z$. On the other hand, notation $xb^k$ is often used as a variant of "scientific" notation of a positive real number in base $b$, with an "integer part"  ( or "first digit") ranging from 1 to $b$-1. In particular, in old, good pre-computer times, one had to use this notation before appealing to  logarithmic tables for base $b$ =10.

A one-to-one map $x \in [1,b) \,|\!\longrightarrow \dfrac{\ln x}{\ln b} \in [0,1)$ realizes an isomorphism between  $Q_b$ and abelian group under addition (mod 1) on $[0,1)$ then easily identifiable with 1-D torus $T$ (and hence circle group and SO(2)).

Let us consider continuous probability density functions $\rho_b$ (pdf) on $[1,b)$ and related cumulative density functions $F_b$ :

$$\int_1^b dx \rho_b(x) = 1 , \qquad F_b(x) = \int_1^x d\tilde{x} \rho_b(\tilde{x}) . \tag{1}$$

There is one special pdf, which will be denoted by $\rho_b^{NB}$ (to honor *both*, S. Newcomb and F. Benford):

$$\rho_b^{NB}(x) = \frac{1}{\ln b} \frac{1}{x} , \text{ with } F_b^{NB}(x) = \frac{\ln x}{\ln b} . \tag{2}$$

For this density, the probability that the first digit in the "scientific notation" using base b is equal 1 will be given by $F_b^{NB}(2)$ (a trivial result for base 2, but leading to $\log_{10} 2 \approx 0.301$ for base 10. The original "Benford law" for occurrence of other digits 2,3,…,9 as "leading ones" follows: Prob( digit 2)= $F_b^{NB}(3) - F_b^{NB}(3) = \log_{10} \frac{3}{2}$ etc.

Clearly, cdf given by (2) also allows to find probabilities of any sequence of digits in arbitrary base as related to the differences of logarithms at that base ( so-called "base invariance" of such cdf).

As $F_b^{NB}$ is just the discussed above isomorphism then, in terms of probability densities, $\rho_b^{NB}$ corresponds to uniform density on [0,1) – a fact noticed early in the history of Benford's law.

Another well-known important property of Newcomb-Benford distribution is its "scale invariance". In this, rather elementary, discussion we will only demonstrate it for intervals $[x_1, x_2] \subset [1, b)$. For any such interval its "Newcomb-Benford measure" is

$$\mu_b^{NB}([x_1, x_2]) = F_b^{NB}(x_2) - F_b^{NB}(x_1) = \frac{1}{\ln b} \ln \frac{x_2}{x_1} . \tag{3}$$

For any $\lambda \in [1, b)$ let us denote by $\lambda[x_1, x_2] (\mathrm{mod}\, b)$ images of intervals $[x_1, x_2]$ under multiplication by $\lambda$ ( mod b ). Simple check of three possible cases:
$\lambda[x_1, x_2](\mathrm{mod}\, b) = [\lambda x_1, \lambda x_2]$, $\lambda[x_1, x_2](\mathrm{mod}\, b) = [\lambda x_1, b) \bigcup [1, x_2 / b)$ and

$\lambda[x_1, x_2](\mathrm{mod}\, b) = [\frac{x_1}{b}, \frac{x_2}{b}]$ leads to

$$\mu_b^{NB}(\lambda[x_1, x_2](\mathrm{mod}\, b)) = \mu_b^{NB}([x_1, x_2]) , \tag{4}$$

i.e. the invariance of the corresponding interval measure.

At this point one should stress that the demonstrated invariance is just a rather elementary example of the invariance of Haar measures on locally compact abelian groups. This fact in the context of the Benford's law was indicated, e.g., in Terence Tao informal blog [3]. General theory of Haar measures gives us also for granted the uniqueness of such invariant measures. A collection of significant results in this direction can be found in [4].

A scheme of working with quotient groups also encourages consideration of "wrapped distributions" when probability distributions on $\mathbb{R}^+$ are "condensed" on the quotient group.

In our case one may start with *non-uniform* partition of $\mathbb{R}^+ = \bigcup_{k \in Z} [b^k, b^{k+1})$ and some cumulative distribution function on $\mathbb{R}^+$. Let us consider a formal series

$$\sum_{k \in Z} \left( F(xb^k) - F(b^{k-1}) \right) \quad . \tag{5}$$

If series (5) is uniformly convergent on $[1, b)$ to some $F_b$, and $F$ itself has an associated density $\rho$, then one may consider a "condensed" or "lumped" density on $[1, b)$:

$$\rho_b(x) = \frac{dF_b}{dx} = \sum_{k \in Z} b^k \rho(xb^k) \quad . \tag{6}$$

To illustrate working of "condensed" densities let us consider the case of log-normal distribution with density

$$\rho_l(y) = \frac{1}{s\sqrt{2\pi}\, y} e^{-(\ln y - M)^2 / 2s^2} \quad , \; y \in R^+ \quad , \tag{7}$$

and "location" and "scale" parameters $M, s$. The corresponding wrapped log-normal distribution is

$$\rho_b^{(l)}(x) = \frac{1}{s\sqrt{2\pi}\, x} \sum_{k \in Z} \frac{1}{b^k} b^k e^{-(\ln x + k \ln b - M)^2 / 2s^2} = \frac{1}{x} \frac{1}{s\sqrt{2\pi}} \sum_{k \in Z} e^{-(\ln x + k \ln b - M)^2 / 2s^2} \quad . \tag{8a}$$

One may try approximate $\rho_b^{(l)}$ by the leading term of the Euler-McLaurin formula, i.e. simply replace the sum in Eq. (8) by an integral with integration variable $k \ln b$. This results in

$$\rho_b^{(l)}(x) \cong \frac{1}{x} \frac{1}{\ln b} \frac{1}{s\sqrt{2\pi}} \int_{-\infty}^{\infty} du \; e^{-(u + \ln x - M)^2 / 2s^2} = \frac{1}{x} \frac{1}{\ln b} = \rho_b^{NB}(x) \quad . \tag{8b}$$

To summarize, Benford's Law will be approximately observed for real positive numbers written in "scientific notation", if these numbers are log-normal distributed with some arbitrary values of the parameters $M, s$. This observation can be immediately extended onto data governed by mixtures of such log-normal distributions.

Log-normal distributions are often observed in "real world"- again with very voluminous body of work on theory and applications [5]. The main reason for that is the frequent occurrence of processes described by products of many independent random variables and related variant of the central limit theorem.

A comment on entropies:

For a density $\rho$ on some set $X$ its entropy $H[\rho]$ is (proportional to, defined as)

$-\int\limits_{X} dx\rho(x)\ln\rho(x)$. For the Newcomb-Benford distribution

$$H[\rho_b^{NB}] = \ln(\ln b) + <\ln x>_{\rho_b^{NB}} = \ln(\ln b) + \frac{1}{2}\ln b \qquad (9)$$

where $<...>_{\rho_b^{NB}}$ is an average calculated with $\rho_b^{NB}$.

For any two pdf 's $\rho, \rho_{ref}$ with common support $X$, a well-known inequality [6 ] reads

$$H[\rho] \leq -\int\limits_{X} dx\rho(x)\ln\rho_{ref}(x) = -<\ln\rho_{ref}(x)>_{\rho}. \qquad (10)$$

Selecting $\rho_b^{NB}$ as reference density leads to

$$H[\rho] \leq \ln(\ln b) + <\ln x>_{\rho} \quad . \qquad (11)$$

Taking into account Eq. (9) one sees that Newcomb-Benford distribution has maximum entropy among all distributions with

$$<\ln x>_{\rho} \leq \frac{1}{2}\ln b. \qquad (12)$$

Provided that our replacement of the infinite series by an integral (Eq. 8a and 8b) is an estimate from below, one can show by the same approximate summation method that the average of ln x calculated with any wrapped log-normal density satisfies the condition of Eq. (12).

General, deep results on entropy of Haar probability measures can be found in [7].

Benford's law is also observed for some sequences generated by deterministic dynamical schemes [8], [2]. One may anticipate that such phenomena could also be related to evolution of probability densities for such systems in the spirit of [6].

(*) Jacek.Kowalski@unt.edu

Wrapped distributions